\documentclass[reqno]{amsproc}
\usepackage[latin1]{inputenc}
\usepackage{amssymb}
\usepackage{amsmath}
\usepackage{latexsym}%
\usepackage[frame,ps,matrix,arrow,curve,rotate,all,2cell,tips]{xy}
\usepackage{cite}
\usepackage{pgf,tikz}
\usetikzlibrary{arrows}
\usepackage{amssymb}
\usepackage{amsfonts}
\usepackage{amscd}
\usepackage{xcolor}
\usepackage{hyperref}
\usepackage{amstext,amsmath,amssymb,amsfonts}
\newtheorem{theorem}{Theorem}

\newtheorem{definition}[theorem]{Definition}
\newtheorem{example}[theorem]{Example}

\newtheorem{lemma}[theorem]{Lemma}

\newtheorem{proposition}[theorem]{Proposition}
\newtheorem{remark}[theorem]{Remark}

\newcommand{\K}{\mathbb {K}}

\newcommand{\A}{\mathcal{A}}

\newcommand{\G}{\mathcal{G}}

\newcommand{\B}{\mathcal{B}}
\newcommand{\bb}{\mathfrak{B}}

\newcommand{\beq}{\begin{eqnarray}}
\newcommand{\eeq}{\end{eqnarray}}
\newcommand{\beqs}{\begin{eqnarray*}}
\newcommand{\eeqs}{\end{eqnarray*}}
\newcommand{\bpro}{\begin{pro}}
\newcommand{\epro}{\end{pro}}
\newcommand{\blem}{\begin{lem}}
\newcommand{\elem}{\end{lem}}
\newcommand{\bdfn}{\begin{dfn}}
\newcommand{\edfn}{\end{dfn}}
\newcommand{\bcor}{\begin{cor}}
\newcommand{\ecor}{\end{cor}}
\newcommand{\bthm}{\begin{thm}}
\newcommand{\ethm}{\end{thm}}
\newcommand{\bex}{\begin{ex}}
\newcommand{\eex}{\end{ex}}
\newcommand{\brmk}{\begin{rmk}}
\newcommand{\ermk}{\end{rmk}}
\newcommand{\bpr}{\begin{pr}}
\newcommand{\epr}{\end{pr}}
\newcommand{\benum}{\begin{enumerate}} 
\newcommand{\eenum}{\end{enumerate}}
\newcommand{\bitem}{\begin{itemize}}
\newcommand{\eitem}{\end{itemize}}

\newcommand{\cqfd}{\hfill{\square}}
\chardef\bslash=`\\
\numberwithin{equation}{section}
\numberwithin{table}{section}
\numberwithin{theorem}{section}
\DeclareMathOperator{\id}{id}

\title{$q$-generalized (anti -) flexible  algebras and  bialgebras}
\author{Mahouton Norbert Hounkonnou$^\ast$}
\address[$\ast$]{University of Abomey-Calavi,
International Chair in Mathematical Physics and Applications,
ICMPA-UNESCO Chair, 072 BP 50, Cotonou, Rep. of Benin}
\email{norbert.hounkonnou@cipma.uac.bj, with copy to hounkonnou@yahoo.fr}
\author{Mafoya Landry Dassoundo$^\dagger$}
\address[$\dagger$]{
University of Abomey-Calavi,
International Chair in Mathematical Physics and Applications, ICMPA-UNESCO Chair, 
072 BP 50, Cotonou, Rep. of Benin
}
\email{mafoya.dassoundo@cipma.uac.bj} 
\begin{document}
\begin{abstract}
In this work, we  provide a $q$-generalization of  flexible algebras
 and related  bialgebraic structures, including  
center-symmetric  (also called antiflexible) algebras, 
and their bialgebras.  
Their basic properties are derived and discussed. Their connection
with known algebraic structures, previously develped 
in the literature, is established.  A $q$-generalization of 
Myung theorem  is given. Main properties related to   bimodules,
matched pairs and dual bimodules as well as their algebraic
consequences are investigated and analyzed. Finally,
the equivalence between $q$-generalized flexible algebras,
their Manin triple and   bialgebras is established.
\end{abstract}
 \maketitle
 


Secondary 16S80, 16T25.
\\
 \today
\tableofcontents
 \section{Introduction}
 Alternative algebras were  introduced by Zorn\cite{zorn} who  established  their fundamental  identities,
 studied   their nucleus by modifying the characteristic of the field, and investigated their Lie admissibility using the corresponding  Jacoby identity. Furthermore, Zorn  
 derived  their power associativity conditions. Later, Schafer  \cite{schafer1} gave a new formulation of 
 these algebras in terms of  left and right multiplication operators, and  in  terms of division algebras of degree two. He also provided 
  the isotopes of these algebras.  Santilli   \cite{santilli} introduced
 Lie admissible algebras and  gave their basic properties. He extended  his study  to  mutation algebras,  examined their relation  to associative algebras, Lie algebras, Jordan and special Jordan algebras,  and established  the passage from   one type of algebra  to another by using 
 a hexahedron with oriented edges. Radicals of flexible Lie admissible algebras were  introduced, and  some of their properties were established
 and discussed in \cite{benkart}. Classes of flexible Lie admissible algebras were also  investigated and discussed in \cite{myung}.   Albert   in \cite{albert} elaborated fundamental concepts, and studied the isotropy of nonassociative algebras. Simple and semi-simple algebras,  and  their characterization  from nonassociative algebraic structures   were developed and discussed in \cite{albert1}.  For more details,  see   a self-contained book by Shafer \cite{schafer} addressing a nice compilation of basic properties of  nonassociative algebras.
 
 Contrarily to Lie algebras, and except for some classification  based on the characteristics of  
 closed fields (see \cite{bhandari} and references therein), a full  classification of  nonassociative  algebras still remains  a tremendous task.
 Some  interesting 
 properties and  algebraic identities of antiflexible structures were investigated and discussed in \cite{celik}. The properties of simple, semisimple  and nearly semisimple 
 antiflexible algebras were also  derived and analyzed in \cite{rodabaugh, rodabaugh1, rodabaugh2}. 
 
 Among the nonassociative
 algebras,  the alternative algebras,
  with an  associator preserving certain symmetry by exchanging its elements,  play a central role  in both mathematics and physics as they possess interesting properties.
 A nice repertory of  their applications in physics, including   gauge theory and  Yang-Mills gauge theory  formulation from nonassociative algebras can be found in \cite{okubo}, (and references therein).  
 A 
 theory of nuclear boson-expansion for odd-fermion system in the context of nonassociative algebras was also examined  in\cite{okubo1}. 
  A  study on   quark structure and  octonion algebras was performed in
 \cite{gunaydin}.
 Further,  a  generalization of the classical Hamiltonian dynamics to a three-dimensional phase space,  generating equations of  motion with two 
 Hamiltonians and three canonical variables,  was performed with  an analog of Poisson bracket  realized by means of the  associator of nonassociative algebras\cite{nambu}.

 Besides, flexible algebras were also investigated in terms of degree of algebras\cite{erwin0}. Other characterizations and applications of nonassociative algebras can be found in \cite{mayne}, \cite{kosier} and \cite{hkn} (and references therein).

 Similarly to algebraic properties of quantum groups developed by Drinfeld~\cite{drinfeld}, 
 some  nonassociative algebras possess interesting identities with applications 
  in physics, and  generate the so-called
 associative or classical Yang-Baxter equations 
 \cite{aguiar, jimbo, drinfeld1},  (and references therein).  Furthermore, the bialgebras constructed from 
 Jordan algebras\cite{Zhelyabin} are related
 to the Lie bialgebras. The left-symmetric
 algebras, also called pre-Lie algebras\cite{bai1},  are known as Lie admissible algebras, and  admit the left multiplication operator
 as a representation. They can also be used to produce symplectic Lie algebras, while their coboundary bialgebras lead to the identity
 known as $S$-equation, and generate para-K\"ahler Lie algebras.
 The case of associative algebras also furnished remarkable properties investigated  by Aguiar\cite{aguiar} and 
  Bai\cite{bai}. The center-symmetric algebras  studied in  \cite{hkn} 
  are also  Lie admissible algebras.

 The present work addresses a $q$-generalization of  flexible algebras and related  bialgebraic structures, including  center-symmetric  (also called antiflexible) algebras, and their bialgebras.  Their basic properties are derived and discussed. Their connection with known algebraic structures existing in the literature is established.  A $q$-generalization of Myung theorem  is given. Main properties related to   bimodules, matched pairs and dual bimodules,  and their algebraic consequences are investigated and analyzed. Finally,   the Manin triple of  $q$-generalized flexible algebras, and its link to $q$-generalized flexible bialgebras are buit together with   the equivalence with the matched pair of  $q$-generalized flexible algebras.
 
 \section{Basic properties of a $q$-generalized flexible algebra}
 In this section,   a $q$-generalization   of algebras encompassing  flexible, anti-flexible and associative  algebras is provided. New classes of algebras are induced. Their relevant properties  and  link with  known algebras are derived. Jordan identity and  Lie admissibility condition  are also established.
 \begin{definition}
 Let $\A(q),$ where $q\in \K,$  be a finite dimensional vector space.  The couple $(\A(q), \cdot )$ 
 is called a $q$-generalized flexible algebra if,    
 for all
 
  $ x, y, z\in \A(q),$ the following relation is satisfied:
 $$ (x, y, z)=q(z, y, x),$$ 
 or, equivalently,
 \begin{eqnarray*}
 \mu\circ(\mu \otimes\id)+q(\mu\circ \tau)\circ ((\mu\circ \tau)\otimes\id)=
 \mu\circ(\id \otimes\mu)+
 q(\mu\circ \tau)\circ (\id \otimes (\mu\circ \tau))
 \end{eqnarray*} 
 where $(x, y, z):=(x\cdot y)\cdot z-x\cdot(y\cdot z)$ is  the associator of the bilinear product "$\cdot$" on $\A(q);$ $\mu$ is defined by $\mu(x, y)=x\cdot y;$ $\id$ is the identity map on $\A(q);$ and $\tau$ stands for the exchange map on $\A(q)$ given by $\tau(x\otimes y)=y\otimes x$.
 \end{definition}
 This can be described  by the following commutative diagram:
 \begin{eqnarray*}
 \label{diagram1}
 \SelectTips{cm}{10}
 \xymatrix{
 \A(q) \otimes \A(q) \otimes \A(q) \ar[rr]^-{\mu\otimes \id-\id\otimes \mu} \ar[d]_-{(\mu\circ \tau)\otimes\id-\id\otimes(\mu\circ\tau)}& & \A(q) \otimes \A(q) \ar[d]^-{\mu}\\
 \A(q)\ar[rr]^-{q\mu\circ\tau} & & \A(q)}
 \end{eqnarray*} 
 
 \begin{remark} We have:
 \begin{itemize}
 \item For $q=0$, the algebra $\A(q)$ is reduced to an associative algebra;
 \item For $q=-1$,   $\A(q)$ becomes a flexible algebra;
 \item For $q=1$,   $\A(q)$ turns to be a center-symmetric algebra  \cite{hkn},  (also called nonflexible algebra \cite{kosier}).
 \end{itemize}
 \end{remark}
 In the sequel, $(\A(q), \cdot)$ denotes a $q$-generalized flexible algebra over $\K$. 
 Besides, for notation simplification, we write $xy$ instead of $x\cdot y$, for $x, y \in \A(q),$ i.e.  the 
  product "$\cdot$" is omitted when there is no confusion.
 \begin{definition}\label{dualsbimodule}
 Suppose $L$ and $R$ be  
 left and right multiplication operators defined on $\A(q)$ as:
 \begin{eqnarray}
 L:
 \begin{array}{lll}
 \A(q) &\longrightarrow \mathfrak{gl}(\A(q)) \\
 x&\longmapsto L_x: 
 \begin{array}{lll}
 \A(q) &\longrightarrow& \A(q) \\
 y &\longmapsto& L_x(y):=x\cdot y
 \end{array}
 \end{array}
 \end{eqnarray}
 \begin{eqnarray}
 R:
 \begin{array}{lll}
 \A(q) &\rightarrow \mathfrak{gl}(\A(q)) \\
 x&\mapsto R_x: 
 \begin{array}{lll}
 \A(q) &\rightarrow& \A(q) \\
 y &\mapsto& R_x(y):=y\cdot x.
 \end{array}
 \end{array}
 \end{eqnarray}
 Then,  their associated dual maps are  defined as follows:
 \begin{eqnarray}
 L^*:
 \begin{array}{lll}
 \A(q) &\rightarrow
 \mathfrak{gl}(\A(q)^*) \\
 x&\mapsto L^*_x: 
 \begin{array}{lll}
 \A(q)^* &\rightarrow& \A(q)^*
  \\
 a &\mapsto& L^*_x(a):
 \begin{array}{lll}
 \A(q) &\rightarrow& \K \cr
 x &\mapsto& 
 \begin{array}{ll}
 <L^*_x(a), y>=\cr 
 \pm <a, L_x(y)>
 \end{array}
 \end{array}
 \end{array}
 \end{array}
 \end{eqnarray}
 \begin{eqnarray}
 R^*:
 \begin{array}{lll}
 \A(q) &\rightarrow
 \mathfrak{gl}(\A(q)^*) \\
 x&\mapsto R^*_x: 
 \begin{array}{lll}
 \A(q)^* &\rightarrow& \A(q)^*
  \\
 a &\mapsto& R^*_x(a):
 \begin{array}{lll}
 \A(q) &\rightarrow& \K \\
 x &\mapsto& 
 \begin{array}{lll}
 <R^*_x(a), y>=\cr 
 \pm <a, R_x(y)>
 \end{array}
 \end{array}
 \end{array}
 \end{array}
 \end{eqnarray}
 \end{definition}
 \begin{proposition}\label{prop_relation1}
 Let 
 $L$ and $R$ be the above defined
 left and right 
 multiplication operators. The  following relations are satisfied for all $x, y \in \A(q):$
 \begin{eqnarray}\label{rel1}
  L_{xy}-L_xL_y&=&q(R_xR_y-R_{yx}), 
 \end{eqnarray}
 \begin{eqnarray}\label{rel2}
 [R_x, L_y]&=&q[R_y, L_x],
 \end{eqnarray}
 \begin{eqnarray}\label{rel3}
  R_xR_y-R_{yx}&=&q(L_{xy}-L_xL_y).
 \end{eqnarray}
 \end{proposition}
 {Proof:}
 
 Let $\A(q)$ be a  $q$-generalized flexible algebra over the field $\K.$
 For all  $x, y, z \in \A(q),$ the proof follows from the equivalences: 
 \begin{eqnarray*}
 (x, y, z)=q(z, y,x) &\Longleftrightarrow& (xy)z-x(yz)=q(zy)x-qz(yx) \cr
 (x, y, z)=q(z, y,x) &\Longleftrightarrow& (L_{xy}-L_xL_y)(z)=q(R_xR_x-R_{yx})(z) \cr
 (x, y, z)=q(z, y,x) &\Longleftrightarrow& (R_zL_x-L_xR_z)(y)=q(R_xL_z-L_zR_x)(y) \cr
 (x, y, z)=q(z, y,x) &\Longleftrightarrow& [R_z, L_x](y)=q[R_x, L_z](y) \cr
 (x, y, z)=q(z, y,x) &\Longleftrightarrow& (R_zR_y-R_{yz})(x)=q(L_{zy}-L_zL_y)(x)
 \end{eqnarray*}
 .$\cqfd$
 \begin{proposition}\label{subajacent}
  Provided the sub-adjacent   algebra 
  $\G(\A(q)):=(\A(q), [, ]),$  where the bilinear product
 $[, ]$ is the  commutator associated to the product  on $\A(q),$  we have,  for all $ x, y, z \in \A(q)$:
 \begin{eqnarray}\label{eqjas}
 J(x, y, z):= [x, [y, z]]+[y, [z, x]]+[z, [x, y]]=\\(q-1)\{(x, y, z)+(y, z, x)+(z, x, y)\}.\nonumber
 \end{eqnarray}
 \end{proposition}
 {Proof:} It stems from  a straightforward computation.
 $\cqfd$
 \begin{lemma}
 \begin{enumerate}
  \item From the Proposition\ref{subajacent}, 
  for all $x, y, z \in \A(q),$ the relation
  \begin{eqnarray}
  S(x, y, z):= (x, y, z)+(y, z, x)+(z, x, y)=0 
  \end{eqnarray}
  is a sufficient condition for $\A(q)$ to become a Lie admissible algebra, i.e. for $(\A(q), [,])$ to  be a Lie algebra. In the particular case where $q=1$,  we get a center symmetric algebra which is Lie admissible  as developed in \cite{hkn}.
 \item The $q$-generalized algebra $\A(q)$ is Lie admissible if and only if,  for all $  x, y, z \in \A(q),$ we have: $(q-1)S(x, y, z)=0.$ In  particular, any $q$-generalized
 flexible algebra defined on a field $\K_{q-1}$ of characteristic $q-1$ is Lie admissible.
 \end{enumerate}
 \end{lemma}
 {Proof:}
 Let $(\A(q), \cdot )$ be a $q$-generalized flexible algebra.
 \begin{enumerate}
  \item From the relation\eqref{eqjas}, we have $J(x,y,z)=(q-1)S(x, y, z).$ Indeed, $\A(q)$ is Lie admissible 
  if and only if $J(x, y, z)=0$  yielding $S(x, y, z)=0$. 
  \item The Lie admissibility condition  $J(x, y, z)=0$  implies that the relation $(q-1)S(x, y, z)=0$ giving
 $S(x, y, z)=0$, or $q-1=0,$ or, in other words, the field $\K$ on which $\A(q)$ defines a vector space is of characteristic $q-1$.$\cqfd$
 \end{enumerate}
 \begin{proposition}
 The following relation is satisfied for all $x, y, z \in \A(q):$
 \begin{eqnarray} \label{rell}
 (L_{xy}-L_xL_y+R_xL_y-L_yR_x+R_yR_x-R_{xy})
 =q(R_xR_y-R_{yx}+R_yL_x\cr -L_xR_y +L_{yx}-L_yL_x),
 \end{eqnarray}
 where $L$ and $R$ are representations of left and right multiplication operators, respectively.
 \end{proposition}
 {Proof:} 
 
 Let us write the associator with the operators $L$ and $R$. 
 
 For all $x, y, z \in \A(q)$,
 \begin{eqnarray*}
 (x, y, z)
 =(xy)z-x(yz) 
 =
 (L_{xy}-L_xL_y)(z) 
 =
 (R_zL_x-L_xR_z)(y)\\
 =(R_zR_y-R_{yz})(x).
 \end{eqnarray*}
 It  follows that:
 \begin{eqnarray*}
 (x, y,z)+(y, z, x)+(z, x, y)&=&(L_{xy}-L_xL_y+R_xL_y-L_yR_x+R_yR_x\\
 &&-R_{xy})(z)\cr
 &=&q(z, y, x)+q(x, z, y)+ 
  q(y, x, z)\cr &=& q(R_xR_y-R_{yx}+R_yL_x
  -L_xR_y+L_{yx}\cr&&- L_yL_x)(z).
 \end{eqnarray*}
 Therefore,  for all $ x, y, z \in \A(q)$:\\
 $L_{xy}-L_xL_y+R_xL_y-L_yR_x+R_yR_x-R_{xy}=q(R_xR_y-R_{yx}+R_yL_x-L_xR_y \\
 +L_{yx}-L_yL_x).$
 $\cqfd$
 \begin{remark}
 The result \eqref{rell} can also be derived  from  the Proposition~\ref{prop_relation1} by summing the relations
 \eqref{rel1}, \eqref{rel2} and \eqref{rel3}.
 \end{remark}
 \begin{theorem}\label{lieadmissibility}
   The following relation is satisfied: $\forall x, y, z \in \A(q),$
  \begin{eqnarray}\label{generalized_jacobi_identity}
 [xy-qyx, z]+[yz-qzy, x]+[zx-qxz, y]=0,
  \end{eqnarray}
  where the bilinear product $[ , ]$ is the commutator associated to the product "$\cdot$" defined on $\A(q)$.
 \end{theorem}
 {Proof:} By a direct computation.
 $\cqfd$
 \begin{remark}
 By setting the parameter  $q=1,$ 
 we get the Jacoby identity from the relation \eqref{generalized_jacobi_identity} indicating  that
  the underlying algebra is Lie admissible  as shown in \cite{hkn}.
 \end{remark}
 We are therefore in right to set the following:
 \begin{definition} Setting $\G(\A(q)):=( \A(q), [,]),$ where $\A(q)$  is the underlying vector space associated to  a
 $q$-generalized flexible algebra
  $(\A(q), \cdot),$
 the equation
 \eqref{generalized_jacobi_identity} defines
 a $q$-generalized Jacobi identity. 
 \end{definition}
 \begin{theorem}\label{generalized_myung_theorem1}
 For a $q$-generalized flexible algebra $ \A(q),$ the following propositions are equivalent: For all $x, y, z \in \A(q),$
 \begin{enumerate}
 \item
 \begin{eqnarray}\label{eq_qderivation}
 [z, xy]=\{z, x\}_qy-x\{y, z\}_q,
 \end{eqnarray}
 where $\{x, y\}_q:=xy+qyx;$
 \item
 \begin{eqnarray}\label{eq_statderivation}
 [z, x\ast_qy]=[z, x]\ast_q y+x\ast_q[z, y],
 \end{eqnarray}
 where $x\ast_q y=\frac{1}{2}(xy-qyx);$
 \item
  $\A(q)$ is a  Lie-admissible $q$-generalized flexible algebra,  i.e.
 \begin{eqnarray}
 [[x, y], z]+[[y, z], x]+[[z, x], y]=0.
 \end{eqnarray}
 \end{enumerate}
 \end{theorem}
 \begin{remark}
  From Theorem~\ref{generalized_myung_theorem1},  we observe that:
 \begin{enumerate}
  \item For $q=-1$,  the equation \eqref{eq_qderivation} turns out to be the derivation for the commutator of a flexible algebra as postulated by 
  the well known Myung Theorem \cite{myung}, \cite{okubo} (and references therein). 
  \item For $q=0$,  the equation \eqref{eq_qderivation} becomes an evidence by using the associativity.
  \item For $q=1$,  the equation \eqref{eq_qderivation} leads to the relation $(x, y, z)+(y, z, x)+(z, x, y)=0,$ which 
 characterizes  an  anti-flexible algebra structure, see  \cite{hkn} (and references therein)
   \item For $q=1$, the equation \eqref{eq_statderivation} is equivalent to the Jacoby identity in a field of a characteristic
   $0$, what is the case for a center-symmetric (also called anti-flexible) algebra;
   \item For $q=0$, the equation \eqref{eq_statderivation} describes the derivation property of    the commutator (or Lie bracket)  of a Lie algebra induced by an
   associative algebra;
   \item For $q=-1$,  the flexibility condition \eqref{eq_statderivation}  defines
   the derivation property of the  Jordan product given as $x\cdot y=\frac{1}{2}(xy+yx)$, see \cite{okubo}.
  \end{enumerate}
 \end{remark}
 \section{Bimodules and matched pairs of $q$-generalized flexible algebras}
 \begin{definition}
 The triple $(l, r, V),$ where $V$ is a finite dimensional vector space and 
 $l, r: \A(q) \rightarrow \mathfrak{gl}(V)$ are two linear maps satisfying the following relations for all $ x, y \in \A(q):$
 \begin{eqnarray}\label{eq_qbimodule1}
 l_{xy}-l_xl_y=q(r_xr_y-r_{yx}),
 \end{eqnarray}
 \begin{eqnarray}\label{eq_qbimodule2}
 [r_x, l_y]=q[r_y, l_x],
 \end{eqnarray}
 \begin{eqnarray}\label{eq_qbimodule3}
 r_xr_y-r_{yx}=q(l_{xy}-l_xl_y),
 \end{eqnarray}
   is called a bimodule of $\A(q),$ also simply denoted by $(l, r).$
 \end{definition}
 
 \begin{proposition}
 Let $l, r: \A(q) \rightarrow \mathfrak{gl}(V)$ be two linear maps. 
 The couple $(l, r)$ is a bimodule of the $q$-generalized flexible algebra $\A(q)$ if and only if  there exists a $q$-generalized flexible 
 algebra structure "$\ast$" on the semi-direct vector space $\A(q)\oplus V$ given by\\ 
 $(x+u)\ast(y+v):=x\cdot y+l_xv+r_yu, \forall x, y \in \A(q)$,  and  $\forall u, v \in V.$ 
 \end{proposition}
 We denote such a   $q$-generalized flexible 
 algebra structure "$\ast$" on the semi-direct vector space $\A(q)\oplus V$ by
 $(\A(q)\oplus V, \ast)$ or simply $\A(q)\ltimes V.$
 
 {Proof:}
 
 Let $x, y, z \in \A(q),$ where $\A(q)$ is a generalized flexible algebra, and $u, v, w \in V$, where $V$
 is a finite dimensional vector space. Using the bilinear product defined, for all $x, y\in \A(q),$ and all $u,v\in V$ by:\\
   $(x+u)\ast(y+v):=x\cdot y+l_xv+r_yu$, where $l,r: \A\rightarrow \mathfrak{gl}(V)$ are linear maps,  the associator of the bilinear product $\ast$  can be written  as:
 \begin{eqnarray*}
 (x+u, y+v, z+w)&=&((x+u)\ast(y+v))\ast(z+w)\cr 
 &&-(x+u)\ast((y+v)\ast(z+w))\cr
 &=&(x\cdot y+l_xv+r_yu)\ast(z+w)\cr 
 &&- (x+u)\ast(y\cdot z+l_yw+r_zv)\cr
 &=&(x\cdot y)\cdot z +l_{x\cdot y}w+r_z(l_xv)+r_z(r_yu)-x\cdot (y\cdot z)\cr && -l_x(l_yw)
 +l_x(r_zv)-r_{y\cdot z}u\cr 
 (x+u, y+v, z+w)&=& (x, y, z)+(l_{x\cdot y}-l_xl_y)w+(z_zl_x-l_xr_z)v\cr &&+(r_zr_y-r_{y\cdot z})u.
 \end{eqnarray*}
 Then, we have, $\forall x, y, z\in \A(q)$ and $\forall u, v,w \in V:$
 \begin{eqnarray}\label{eq_qbj1}
 (x+u, y+v, z+w)=(x, y, z)+(l_{x\cdot y}-l_xl_y)w+(z_zl_x-l_xr_z)v\cr 
 +(r_zr_y-r_{y\cdot z})u.
 \end{eqnarray}
 Besides,
 \begin{eqnarray}\label{eq_qbj2}
 q(z+w, y+v, x+u)=q(z, y, x)+q(l_{z\cdot y}-l_zl_y)u+q(z_xl_z-l_zr_x)v\cr 
 +q(r_xr_y-r_{y\cdot x})w.
 \end{eqnarray}
 By setting $(x+u, y+v, z+w)=q(z+w, y+v, x+u)$, $\forall x, y, z\in \A(q)$ and $\forall u, v,w \in V$,   i.e. $(\A(q)\oplus V, \ast)$ is
 a $q$-generalized flexible algebra, we get from the right hand side of the   equations \eqref{eq_qbj1} and \eqref{eq_qbj2}:
 \begin{eqnarray*}
 (x, y, z)&=&q(z, y, x),\\
 (l_{x\cdot y}-l_xl_y)&=&q(r_xr_y-r_{y\cdot x}),\\
 (z_zl_x-l_xr_z)&=&q(z_xl_z-l_zr_x),\\
 (r_zr_y-r_{y\cdot z})&=&q(l_{z\cdot y}-l_zl_y),
 \end{eqnarray*}
 which are equivalent to the relations \eqref{eq_qbimodule1}, \eqref{eq_qbimodule2}, \eqref{eq_qbimodule3} defining the bimodule of a $q$-generalized flexible algebra. 
 $\cqfd$
 \begin{example}
 According to the Proposition\ref{prop_relation1},   $(L, R),$ where $R$ and $L$ are the representations of the  right and left multiplication operators, respectively,  is a bimodule of a $q$-generalized flexible
 algebra $\A(q).$ Indeed, $L$ and $R$ satisfy
 the equations \eqref{eq_qbimodule1}, \eqref{eq_qbimodule2} and \eqref{eq_qbimodule3}.
 \end{example}
 \begin{theorem}\label{theo_01}
  Let $(\A(q), \cdot)$ and $(\B(q), \ast)$ be two $q$-generalized flexible algebras. Suppose  there  exist linear maps
  $l_{\A}, r_{\A}: \A(q)\rightarrow \mathfrak{gl}(\B(q))$ and $l_{\B}, r_{\B}: \B(q)\rightarrow \mathfrak{gl}(\A(q))$ 
  satisfying the following relations:
  \begin{eqnarray}\label{eq_qmatch1}
  (l_{\B}(a)x)\cdot y+l_{\B}(r_{\A}(x)a)y-l_{\B}(a)(x\cdot y)=q(r_{\B}(a)(y\cdot x)-y\cdot(r_{\B}(a)x)\cr -r_{\B}(l_{\A}(x)a)y),
  \end{eqnarray}
  \begin{eqnarray}\label{eq_qmatch2}
 r_{\B}(a)(x\cdot y)-x\cdot(r_{\B}(a)y)-r_{\B}(l_{\A}(y)a)x=q((l_{\B}(a)y)\cdot x+l_{\B}(r_{\A}(y)a)x\cr 
 -l_{\B}(a)(y\cdot x)),
  \end{eqnarray}
  \begin{eqnarray}\label{eq_qmatch3}
  &&(r_{\B}(a)x)\cdot y+l_{\B}(l_{\A}(x)a)y-x\cdot(l_{\B}(a)y)-r_{\B}(r_{\A}(y)a)x=
  q(((r_{\B}(a)y)\cdot x\cr 
  && +l_{\B}(l_{\A}(y)a)x- y\cdot (l_{\B}(a)x)-r_{\B}(r_{\A}(x)a)y),
  \end{eqnarray}
  \begin{eqnarray}\label{eq_qmatch4}
 &&(l_{\A}(x)a)\ast b+l_{\A}(r_{\B}(a)x)b-l_{\A}(x)(a\ast b)=q(r_{\A}(x)(b\ast a) 
 -b\ast (r_{\A}(x)a)\cr 
 &&- 
 r_{\A}(l_{\B}(a)x)b)),
  \end{eqnarray}
   \begin{eqnarray}\label{eq_qmatch5}
 && r_{\A}(x)(a\ast b)-a\ast (r_{\A}(x)b)-r_{\A}(l_{\B}(b)x)a=q((l_{\A}(x)b)\ast a+l_{\A}(r_{\B}(b)x)a\cr
 && 
 -l_{\A}(x)(b\ast a)),
  \end{eqnarray}
  \begin{eqnarray}\label{eq_qmatch6}
  &&(r_{\A}(x)a)\ast b+l_{\A}(l_{\B}(a)x)b-a\ast(l_{\A}(x)b)-r_{\A}(r_{\B}(b)x)a=
  q((r_{\A}(x)b)\ast a\cr 
 &&+l_{\A}(l_{\B}(b)x)a-  b\ast (l_{\A}(x)a)-r_{\A}(r_{\B}(a)x)b)),
  \end{eqnarray}
 $\forall x, y\in \A(q)$ and $\forall a, b\in\B(q).$ It  follows that  there is a $q$-generalized flexible algebra structure "$\star$" on the direct sum of vector spaces
 $\A(q)\oplus\B(q)$ given as: $(x+a)\star(y+b)=(x\cdot y+l_{\B}(a)y+r_{\B}(b)x)+(a\ast b+l_{\A}(x)b+r_{\A}(y)a).$ 
 \end{theorem}
 
 {Proof:}
 
 Let $(\A(q), \cdot )$ and $(\B(q), \ast)$ be two $q$-generalized flexible algebras,
 $l_{\A}, r_{\A}: \A(q) \rightarrow \mathfrak{gl}(\B(q))$ and 
 $l_{\B}, r_{\B}: \B(q) \rightarrow \mathfrak{gl}(\A(q))$ be four linear maps satisfying the relations
 \eqref{eq_qmatch1}, \eqref{eq_qmatch2}, \eqref{eq_qmatch3}, \eqref{eq_qmatch4}, \eqref{eq_qmatch5} and 
 \eqref{eq_qmatch6}. Consider the bilinear product "$\star$"
 defined  on the vector space  
 $\A(q)\oplus\B(q)$ as: $(x+a)\star(y+b)=(x\cdot y+l_{\B}(a)y+r_{\B}(b)x)+(a\ast b+l_{\A}(x)b+r_{\A}(y)a), \forall x, y\in \A(q); a,b \in \B(q).$ We have:
 \begin{eqnarray*}
 (x+a, y+b, z+c)&=&\{(x+a)\star(y+b)\}\star(z+c)\cr 
 &&-(x+a)\star\{(y+b)\star(z+b)\}\cr
 &=& (x, y, z) + (a, b, c) + \{r_{\mathcal{B}}(c)(x \cdot y) + 
 l_{\mathcal{A}}(x \cdot y)c \cr 
 && -x \cdot (r_{\mathcal{B}}(c)y) 
 -l_{\mathcal{A}}(x)(l_{\mathcal{A}}(y)c)
 -r_{\mathcal{B}}(l_{\mathcal{A}}(y)c)x \}  \cr 
 &&
 +\{r_{\mathcal{B}}(c)(r_{\mathcal{B}}(b)x) 
 + l_{\mathcal{A}}(r_{\mathcal{B}}(b)x)c 
 -r_{\mathcal{B}}(b \ast c)x  \cr &&+ (l_{\mathcal{A}}(x)b)\ast c-
 l_{\mathcal{A}}(x) ( b \circ c) \} 
  \cr
 &&+ \{(r_{\mathcal{B}}(b)x) \cdot z 
  l_{\mathcal{B}}(l_{\mathcal{A}}(x)b)z
  + r_{\mathcal{A}}(z)(l_{\mathcal{A}}(x)b)\cr 
 &&
 -x \cdot (l_{\mathcal{B}}(b)z)- 
  r_{\mathcal{B}}(r_{\mathcal{A}}(z)b)x-l_{\mathcal{A}}(x)(r_{\mathcal{A}}(z)b)\} 
 \cr &&
 + \{(l_{\mathcal{B}}(a)y)\cdot z
  + 
 l_{\mathcal{B}}(r_{\mathcal{A}}(y)a)z  
 +r_{\mathcal{A}}(z)(r_{\mathcal{A}}(y)a)
 \cr 
 &&
 -l_{\mathcal{B}}(a)(y \cdot z)- r_{\mathcal{A}}(y \cdot z )a\}
 +
 \{r_{\mathcal{B}}(c)(l_{\mathcal{B}}(a)y) 
 \cr 
 &&+(r_{\mathcal{A}}(y)a) \ast c +
  l_{\mathcal{A}}(l_{\mathcal{B}}(a)y)c-
 +l_{\mathcal{B}}(a)(r_{\mathcal{B}}(c)y)
 \cr &&
 - a \ast (l_{\mathcal{A}}(y)c) 
 -r_{\mathcal{A}}(r_{\mathcal{B}}(c)y)a\} + \{l_{\mathcal{B}}(a \ast b)z 
  \cr&&
 +  r_{\mathcal{A}}(z)(a \ast b)-
 + l_{\mathcal{B}}(a)(l_{\mathcal{B}}(b)z) 
 - a \ast (r_{\mathcal{A}}(z)b) \cr 
 &&-r_{\mathcal{A}}(l_{\mathcal{B}}(b)z)a\}\cr
 &=&(x, y, z)+(x, y, c)+(x, b, z)+(x, b, c)+(a,y,z)\cr 
 &&+(a, b, z)+(a, b, c)
 \end{eqnarray*}
 Then, the $q$-generalized flexibility condition of the bilinear product $\star$ is given as:
 \begin{eqnarray*}
 (x+a, y+b, z+c)=q(z+c, y+b, x+a)&\Longleftrightarrow& 
 \left\lbrace 
 \begin{array}{lll}
 (x, y, z)=q(z, y, x) \\
 (a, b, c)=q(c, b, a)\\
 (x, y, c)=q(c, y, x)\\
 (x, b, z)=q(z, b, x)\\
 (x, b, c)=q(c, b, x)\\
 (a, y, z)=q(z, y, a)\\
 (a, y, c)=q(c, y, a)\\
 (a, b, z)=q(z, b, a)
 \end{array}
 \right.\cr
 &\Longleftrightarrow &
 \left\lbrace
 \begin{array}{lllll}
 (l_{\A}, r_{\A}, \B(q)), \\ (l_{\B}, r_{\B}, \A(q))\\
 \mbox{ are bimodules of } \\\A(q), \B(q),\\ \mbox{ respectively,}\\
 (x, y, a)=q(a, y, x)\\
 (x, a, y)=q(y, a, x)\\
 (x, a, b)=q(b, a, x)\\
 (a, x, y)=q(y, x, a)\\
 (a, x, b)=q(b, x, a)\\
 (a, b, x)=q(x, b, a)
 \end{array}
 \right.\\
 (x+a, y+b, z+c)=q(z+c, y+b, x+a)& \Longleftrightarrow&
 \left\lbrace 
 \begin{array}{lllll}
 (l_{\A}, r_{\A}, \B(q)),\\ (l_{\B}, r_{\B}, \A(q))\\
 \mbox{ are bimodules of } \\ \A(q), \B(q), \\\mbox{ respectively,}
 \mbox{ and the}\\ \mbox{ following } \mbox{ equations:}
 \\
 \eqref{eq_qmatch1}, \eqref{eq_qmatch2}, \\ \eqref{eq_qmatch3}, 
 \eqref{eq_qmatch4}, \eqref{eq_qmatch5}, \eqref{eq_qmatch6}  \\ \mbox{ are satisfied. }
 \end{array}
 \right.
 \end{eqnarray*}
  Therefore, we obtain a $q$-generalized flexible algebra structure given, for all $x, y\in \A(q)$ and all $a, b\in \B(q),$ by:
    $(x+a)\star(y+b)=(x\cdot y+l_{\B}(a)y+r_{\B}(b)x)+(a\ast b+l_{\A}(x)b+r_{\A}(y)a)$ 
 on the direct sum of the underlying vector spaces  $\A(q)$ and  $\B(q)$ of the bimodules  ($l_{\A}, r_{\A}, \B(q)$) and ($l_{\B}, r_{\B}, \A(q)$)
  of the associated $q$-generalized flexible algebras  $\A(q)$ and $\B(q).$
 $\cqfd$
 
  In this case, 
 the obtained $q$-generalized flexible algebra $(\A(q)\oplus\B(q), \star)$ is denoted by 
 $\A(q)\bowtie^{l_{\A}, r_{\A}}_{l_{\B}, r_{\B}}\B(q),$ or simply $\A(q)\bowtie\B(q).$ 
 \begin{definition}
  The sixtuple 
 $(\A(q), \B(q), l_{\A}, r_{\A}, l_{\B}, r_{\B})$ satisfying the conditions of Theorem~\ref{theo_01} 
  is called matched pair 
  of the generalized flexible algebras $\A(q)$ and $\B(q).$
 \end{definition}

 \begin{remark}
  Theorem~\ref{theo_01} is a $q$-generalization of   main theorems, well known in the literature.  Indeed,
 \begin{itemize}
  \item For $q=0$,  Theorem~\ref{theo_01} is  exactly the fundamental theorem for the  matched pair of associative algebras. See \cite{bai} and references therein.
  \item For $q=1$,  Theorem~\ref{theo_01} is reduced to the fundamental theorem for the matched pair of  center-symmetric algebras formulated in \cite{hkn}.
  \item For $q=-1$,  Theorem~\ref{theo_01} becomes  the fundamental theorem for the matched pair of flexible algebras.
 \end{itemize}
 \end{remark}
 \section{Basic properties of the $q$-generalized flexible algebras}
 In this section, we construct  and discuss the  basic definitions and main properties of the $q$-generalized flexible algebras.
 \begin{definition}
 Let $l, r: \A(q)\rightarrow \mathfrak{gl}(V)$ be the two above mentioned  linear maps.
 Their dual maps are  defined as:
 \begin{eqnarray}
 l^*:
 \begin{array}{lll}
 \A(q) &\rightarrow
 \mathfrak{gl}(V^*) \\
 x&\mapsto l^*_x: 
 \begin{array}{lll}
 V^* &\rightarrow& V^*
  \\
 v^* &\mapsto& l^*_x(v^*):
 \begin{array}{lll}
 V &\rightarrow& \K \\
 u &\mapsto&
 \begin{array}{lll}
  <l^*_x(v^*), u>=\\
 <v^*, l_x(u)>,
 \end{array}
 \end{array}
 \end{array}
 \end{array}
 \end{eqnarray}
 \begin{eqnarray}
 r^*:
 \begin{array}{lll}
 \A(q) &\rightarrow
 \mathfrak{gl}(V^*) \\
 x&\mapsto r^*_x: 
 \begin{array}{lll}
 V^* &\rightarrow& V^*
  \\
 v^* &\mapsto& r^*_x(v^*):
 \begin{array}{lll}
 V &\rightarrow& \K \\
 u &\mapsto& 
 \begin{array}{lll}
  <r^*_x(v^*), u>=\\
 <v^*, r_x(u)>,
 \end{array}
 \end{array}
 \end{array}
 \end{array}
 \end{eqnarray}
 where $V^*=Hom(V, \K)$ and $\mathfrak{gl}(V^*)$ is the linear group of $V^*.$
 \end{definition}
 \begin{theorem}\label{prop_dualbimodule}
  For any finite dimensional vector space $V$, suppose 
 $l, r: \A(q) \rightarrow \mathfrak{gl}(V)$ be two linear maps such that $l^*$ and $r^*$ are dual maps of $l$ and $r,$ respectively. The following propositions are equivalent:
 \begin{enumerate}
 \item $(l, r, V)$ is a bimodule of the $q$-generalized flexible algebra $\A(q),$
 \item $(r^*, l^*, V^*)$ is a bimodule of the $q$-generalized flexible algebra $\A(q).$
 \end{enumerate} 
 \end{theorem}
 \begin{remark}
 It is worth noticing  that the dual bimodule does not depend on the parameter  $q$. Clearly, it is the dual bimodule 
 obtained for the center-symmetric algebra in  \cite{hkn}.  Therefore,  the dual bimodule of center-symmetric  
 algebras is  the same as the dual bimodule of  flexible algebras. 
 \end{remark}
 \begin{proposition}
 The quadruple  $(R_{\cdot}^*, L_{\cdot}^*, A(q)^*),$  where $A(q)^*$ is the dual space of $\A(q)$ and given by $A(q)^*=Hom(\A(q), \K)$,  is a bimodule of $\A(q).$
 \end{proposition}
 {Proof:}
 By considering  Definition~\ref{dualsbimodule},  Proposition~\ref{prop_relation1} and  Proposition~\ref{prop_dualbimodule}, 
 we deduce that $(R_{\cdot}^*, L_{\cdot}^*, A(q)^*)$ is a bimodule of $\A(q).$
 $\cqfd$
 \begin{theorem}\label{theo_dumatchedpair}
 Let $(\A(q), \cdot)$ be a $q$-generalized flexible algebra. Suppose that there is a $q$-generalized flexible algebra structure  "$\circ$" on its dual space 
 $\A(q)^*=Hom(\A(q), \K)$. The  sixtuple $(\A(q), \A(q)^*, R_{\cdot}^*, L_{\cdot}^*, R_{\circ}^*, L_{\circ}^*)$ is a matched pair of the $q$-generalized flexible algebras
 $\A(q)$ and  $\A(q)^*$ if and only if the linear maps $ R_{\cdot}^*, L_{\cdot}^*, R_{\circ}^*, L_{\circ}^*$ satisfy the following relations for all
 $x, y\in\A(q)$ and all $a, b\in \A(q)^*$:
 \begin{eqnarray}\label{eq_qdualmatch1}
 &&(R_{\circ}^*(a)x)\cdot y+R_{\circ}^*(L_{\cdot}^*(x)a)y-R_{\circ}^*(a)(x\cdot y)=q(L_{\circ}^*(a)(y\cdot x)\cr 
 && -y\cdot(L_{\circ}^*(a)x)-L_{\circ}^*(R_{\cdot}^*(x)a)y),
 \end{eqnarray}
 \begin{eqnarray}\label{eq_qdualmatch2}
 && L_{\circ}^*(a)(x\cdot y)-x\cdot(L_{\circ}^*(a)y)-L_{\circ}^*(R_{\cdot}^*(y)a)x=q((R_{\circ}^*(a)y)\cdot x
 \cr
 && +R_{\circ}^*(L_{\cdot}^*(y)a)x
  -R_{\circ}^*(a)(y\cdot x)),
 \end{eqnarray}
 \begin{eqnarray}\label{eq_qdualmatch3}
 &&(L_{\circ}^*(a)x)\cdot y+R_{\circ}^*(R_{\cdot}^*(x)a)y-x\cdot(R_{\circ}^*(a)y)-L_{\circ}^*(L_{\cdot}^*(y)a)x=\cr 
 && 
 q((L_{\circ}^*(a)y)\cdot x+  R_{\circ}^*(R_{\cdot}^*(y)a)x- y\cdot (R_{\circ}^*(a)x)-L_{\circ}^*(L_{\cdot}^*(x)a)y).
 \end{eqnarray}
 \end{theorem}
 {Proof:}
 
 Consider a $q$-generalized flexible algebra $(\A(q), \cdot)$ and 
  assume that there is a $q$-generalized flexible algebra
   structure "$\circ$" on its dual
   space $(\A(q))^*$. Using  Definition~\ref{dualsbimodule},
  Proposition~\ref{prop_relation1} and  Proposition~\ref{prop_dualbimodule}, we deduct that $(R_{\cdot}^*, L_{\cdot}^*, \A(q)^*)$ is a
  bimodule of $\A(q)$ and $(R_{\circ}^*, L_{\circ}^*, \A(q))$ is a bimodule of $\A(q)^*$.  It remains to show that the quadruple 
  ($(R_{\cdot}^*, L_{\cdot}^*, R_{\circ}^*, L_{\circ}^*$) satisfies the relations 
  \eqref{eq_qmatch1}, \eqref{eq_qmatch2}, \eqref{eq_qmatch3}, \eqref{eq_qmatch4}, \eqref{eq_qmatch5} and \eqref{eq_qmatch6}.
 By taking the correspondences 
 $l_{\A}\rightarrow R_{\cdot}^*$, 
 $r_{\A}\rightarrow L_{\cdot}^*$,
 $l_{\B}\rightarrow R_{\cdot}^*$ and
 $r_{\A}\rightarrow L_{\cdot}^*$, we  straightforwardly obtain the relations
 \eqref{eq_qmatch1}, \eqref{eq_qmatch2} and 
 \eqref{eq_qmatch3}. To establish now the relations \eqref{eq_qmatch4}, \eqref{eq_qmatch5} and 
 \eqref{eq_qmatch6}, let us compute, $\forall x, y \in \A(q)$ and $\forall a,b\in \B(q):$ 
 \begin{eqnarray*}
 &&< (R_{\circ}^*(a)x)\cdot y+R_{\circ}^*(L_{\cdot}^*(x)a)y-R_{\circ}^*(a)(x\cdot y)-q(L_{\circ}^*(a)(y\cdot x)\\
 &&-y\cdot(L_{\circ}^*(a)x)-L_{\circ}^*(R_{\cdot}^*(x)a)y), b > 
 =< L_{\cdot}(R_{\circ}^*(a)x) y+R_{\circ}^*(L_{\cdot}^*(x)a)y\\
 && -R_{\circ}^*(a)(L_{\cdot}(x)y)- q(L_{\circ}^*(a)(R_{\cdot}(x)y)
 -R_{\cdot}(L_{\circ}^*(a)x)y-L_{\circ}^*(R_{\cdot}^*(x)a)y), b>\\
 &&
 = <y,  L_{\cdot}^*(R_{\circ}^*(a)x)b+ (L_{\cdot}^*(x)a)\circ b-L_{\cdot}^*(x)(b\circ a)-q(
 (R_{\cdot}^*(x)(a\circ b))\\ && -q(R_{\cdot}^*(L_{\circ}^*(a)x)b)
 +q((R_{\cdot}^*(x)a)\circ b)>.
 \end{eqnarray*}
 Then we get  for all $x, y \in \A(q)$ and all $a, b\in \B(q):$
 \begin{eqnarray*}
 &&< (R_{\circ}^*(a)x)\cdot y+R_{\circ}^*(L_{\cdot}^*(x)a)y-R_{\circ}^*(a)(x\cdot y)-q(L_{\circ}^*(a)(y\cdot x)\cr 
 && -y\cdot(L_{\circ}^*(a)x)-
 L_{\circ}^*(R_{\cdot}^*(x)a)y), b> 
 =<y,  L_{\cdot}^*(R_{\circ}^*(a)x)b+(L_{\cdot}^*(x)a)\circ b\cr 
 && -L_{\cdot}^*(x)(b\circ a)-q(
 (R_{\cdot}^*(x)(a\circ b)-
 R_{\cdot}^*(L_{\circ}^*(a)x)b
 -(R_{\cdot}^*(x)a)\circ b
 )>.
 \end{eqnarray*}
 The following holds:
 \begin{eqnarray*}
 (R_{\circ}^*(a)x)\cdot y+R_{\circ}^*(L_{\cdot}^*(x)a)y-R_{\circ}^*(a)(x\cdot y)&=&q(L_{\circ}^*(a)(y\cdot x)
 -y\cdot(L_{\circ}^*(a)x)\cr 
 &&- L_{\circ}^*(R_{\cdot}^*(x)a)y)
 \end{eqnarray*}
 or, equivalently,
 \begin{eqnarray*}
 L_{\cdot}^*(R_{\circ}^*(a)x)b+ b\circ (L_{\cdot}^*(x)a)-L_{\cdot}^*(x)(b\circ a)=q(
 (R_{\cdot}^*(x)(a\circ b)-
 R_{\cdot}^*(L_{\circ}^*(a)x)b\cr 
 - (R_{\cdot}^*(x)a)\circ b)
 \end{eqnarray*}
 which is  exactly the equation \eqref{eq_qmatch5} by taking the correspondences: $l_{\A}\rightarrow R_{\cdot}^*$, 
 $r_{\A}\rightarrow L_{\cdot}^*$,
 $l_{\B}\rightarrow R_{\cdot}^*$,
 $r_{\A}\rightarrow L_{\cdot}^*$, $a \rightarrow b,$ and we have 
 \begin{eqnarray}\label{eq_qdualmatch5}
 && L_{\cdot}^*(x)(a\circ b)-a\circ (L_{\cdot}^*(x)b)-L_{\cdot}^*(R_{\circ}^*(b)x)a=q((R_{\cdot}^*(x)b)\circ a\cr 
 && +R_{\cdot}^*(L_{\circ}^*(b)x)a  -R_{\cdot}^*(x)(b\circ a)).
 \end{eqnarray}
 This  shows that the relation \eqref{eq_qdualmatch1} is equivalent to \eqref{eq_qdualmatch5}. 
 By the same way, we have:
 \begin{eqnarray*}
 &&< L_{\circ}^*(a)(x\cdot y)-x\cdot(L_{\circ}^*(a)y)-L_{\circ}^*(R_{\cdot}^*(y)a)x-q((R_{\circ}^*(a)y)\cdot x\cr 
 && +R_{\circ}^*(L_{\cdot}^*(y)a)x
 -R_{\circ}^*(a)(y\cdot x)), b >=
 <L_{\circ}^*(a)(R_{\cdot}(y)x)- 
 R_{\cdot}(L_{\circ}^*(a)y)x\cr 
 && -L_{\circ}^*(R_{\cdot}^*(y)a)x-q(L_{\cdot}(R_{\circ}^*(a)y)x+
 R_{\circ}^*(L_{\cdot}^*(y)a)x
  -R_{\circ}^*(a)(L_{\cdot}(y)x)), b  >\cr 
  &&
 = <x,R_{\cdot}^*(y)(a\circ b)-
  R_{\cdot}^*(L_{\circ}^*(a)y)b-
 (R_{\cdot}^*(y)a)\circ b
 -q[L_{\cdot}^*(R_{\circ}^*(a)y)b \cr 
 && + b\circ(L_{\cdot}^*(y)a)-L_{\cdot}^*(y)(b \circ a)] >.
 \end{eqnarray*}
 Then, we have that
 \begin{eqnarray*}
 && L_{\circ}^*(a)(x\cdot y)-x\cdot(L_{\circ}^*(a)y)-L_{\circ}^*(R_{\cdot}^*(y)a)x=q((R_{\circ}^*(a)y)\cdot x\cr 
 && +R_{\circ}^*(L_{\cdot}^*(y)a)x-R_{\circ}^*(a)(y\cdot x))
 \end{eqnarray*}
 is equivalent to
 \begin{eqnarray*}
 && R_{\cdot}^*(y)(a\circ b)-R_{\cdot}^*(L_{\circ}^*(a)y)b-
 (R_{\cdot}^*(y)a)\circ b=q(L_{\cdot}^*(R_{\circ}^*(a)y)b\cr 
 && 
 +b\circ(L_{\cdot}^*(y)a)-L_{\cdot}^*(y)(b \circ a)).
 \end{eqnarray*}
 This  exactly gives the equation \eqref{eq_qmatch4}
 by setting the correspondence 
  $l_{\A}\rightarrow R_{\cdot}^*$, 
 $r_{\A}\rightarrow L_{\cdot}^*$,
 $l_{\B}\rightarrow R_{\cdot}^*$,
 $r_{\A}\rightarrow L_{\cdot}^*$, $y \rightarrow x,$ and then we get:
 \begin{eqnarray}\label{eq_qdualmatch4}
 &&(R_{\cdot}^*(x)a)\circ b+R_{\cdot}^*(L_{\circ}^*(a)x)b-R_{\cdot}^*(x)(a\circ b)=q(L_{\cdot}^*(x)(b\circ a)\cr 
 && -b\circ(L_{\cdot}^*(x)a)-  L_{\cdot}^*(R_{\circ}^*(a)x)b).
 \end{eqnarray}
 This  shows that the relation \eqref{eq_qdualmatch2} is equivalent to \eqref{eq_qdualmatch4}: 
 \begin{eqnarray*}
 &&< (L_{\circ}^*(a)x)\cdot y+R_{\circ}^*(R_{\cdot}^*(x)a)y-x\cdot(R_{\circ}^*(a)y)-L_{\circ}^*(L_{\cdot}^*(y)a)x\\
 &&
  -q((L_{\circ}^*(a)y)\cdot x+R_{\circ}^*(R_{\cdot}^*(y)a)x -y\cdot (R_{\circ}^*(a)x)-L_{\circ}^*(L_{\cdot}^*(x)a)y), b  > =\\
  &&
  <y, L_{\cdot}^*(L_{\circ}(a)x)b+b\circ(R_{\cdot}^*(x)a)
  - (L_{\cdot}^*(x)b)\circ a
 -R_{\cdot}^*(R_{\circ}^*(b)x)a\\
 &&
 -q(a\circ (R_{\cdot}^*(x)b) +L_{\cdot}^*(L_{\circ}^*(b)x)a)-
  R_{\cdot}^*(R_{\circ}^*(a)x)b-(L_{\cdot}^*(x)a)\circ b)>.
 \end{eqnarray*}
 Therefore, the following relation
 \begin{eqnarray*}
 (L_{\circ}^*(a)x)\cdot y+R_{\circ}^*(R_{\cdot}^*(x)a)y-x\cdot(R_{\circ}^*(a)y)-L_{\circ}^*(L_{\cdot}^*(y)a)x=
  q((L_{\circ}^*(a)y)\cdot x \\ +R_{\circ}^*(R_{\cdot}^*(y)a)x-
   y\cdot (R_{\circ}^*(a)x)-L_{\circ}^*(L_{\cdot}^*(x)a)y)
 \end{eqnarray*}
 is equivalent to 
 \begin{eqnarray*}
 L_{\cdot}^*(L_{\circ}(a)x)b+b\circ(R_{\cdot}^*(x)a)
  - (L_{\cdot}^*(x)b)\circ a-
 R_{\cdot}^*(R_{\circ}^*(b)x)a-
   q(a\circ (R_{\cdot}^*(x)b) \\
  + L_{\cdot}^*(L_{\circ}^*(b)x)a)-
  R_{\cdot}^*(R_{\circ}^*(a)x)b-(L_{\cdot}^*(x)a)\circ b)
 \end{eqnarray*}
 which is exactly the equation \eqref{eq_qmatch6}
 by setting the correspondences \\
  $l_{\A}\rightarrow R_{\cdot}^*$, 
 $r_{\A}\rightarrow L_{\cdot}^*$,
 $l_{\B}\rightarrow R_{\cdot}^*$,
 $r_{\A}\rightarrow L_{\cdot}^*$, $a \rightarrow b,$ and we have  
 \begin{eqnarray}\label{eq_qdualmatch6}
  &&(L_{\cdot}^*(x)a)\circ b+R_{\cdot}^*(R_{\circ}^*(a)x)b-a
  \circ(R_{\cdot}^*(x)b)-L_{\cdot}^*(L_{\circ}^*(b)x)a=\cr 
  && 
  q((L_{\cdot}^*(x)b)\circ a+  R_{\cdot}^*(R_{\circ}^*(b)x)a)-  b\circ (R_{\cdot}^*(x)a)-L_{\cdot}^*(L_{\circ}^*(a)x)b)).
 \end{eqnarray}
 This  shows that the relation \eqref{eq_qdualmatch3} is equivalent to \eqref{eq_qdualmatch6}. 
 
 Therefore, the sixtuple $(\A(q), (\A(q)^*, R_{\cdot}^*, L_{\cdot}^*, R_{\circ}^*, L_{\circ}^*)$ is a matched pair of the $q$-generalized flexible  algebras 
 $\A(q)$ and $\A(q)^*$ if and only if the linear maps $ R_{\cdot}^*, L_{\cdot}^*, R_{\circ}^*, L_{\circ}^*$ satisfy the equations 
 \eqref{eq_qdualmatch1}, \eqref{eq_qdualmatch2} and \eqref{eq_qdualmatch3}.
 $\cqfd$
 \begin{remark}
 Theorem~\ref{theo_dumatchedpair} encompasses particular results known in the literature, namely:
 \begin{itemize} 
 \item For $q=0,$   Theorem\ref{theo_dumatchedpair} is exactly reduced to the result  obtained by Bai in \cite{bai}, 
 (see also references therein) giving the  construction of  the dual matched pair for the associative algebras.
 \item For $q=1,$  we recover   the theorem  relating 
 the dual matched pair of center-symmetric algebras with the dual matched pair of Lie algebras, investigated in \cite{hkn}.
 \item For $q=-1,$ Theorem\ref{theo_dumatchedpair}  gives  the dual matched pair of  flexible algebras.  This  is a new result,  given in this work for the first time, to our best knowledge of the literature.
 \end{itemize}
 \end{remark}
 \begin{proposition}\label{prop_sum} 
 Assume that there is a $q$-generalized flexible algebra structure 
 "$\circ$" on the dual space $\A(q)^*$. There is a $q$-generalized flexible algebra structure "$\star$" on the vector space $\A(q)\oplus\A(q)^*$ given, 
 for all $x, y\in\A(q)$ and all  $a, b\in \A(q)^*,$  by:
 \begin{eqnarray}\label{eq_qproddual}
 (x+a)\star(y+b)=(x\cdot y+R_{\circ}^*(a)y+L_{\circ}^*(b)x)\cr 
 + (a\circ b+R_{\cdot}^*(x)b+L_{\cdot}^*(y)a),
 \end{eqnarray} 
 if and only if   the sixtuple $(\A(q), \A(q)^*, R_{\cdot}^*, L_{\cdot}^*, R_{\circ}^*, L_{\circ}^*)$ is a matched pair of the $q$-generalized 
 flexible algebras $\A(q)$ and $\A(q)^*$.
 \end{proposition}
 {Proof:}
 
 It is well known that  $\A(q)\oplus\A(q)^*$ is a vector space as a direct sum of  vector spaces. The product "$\star$" is also a  bilinear product by definition. Then, it only remains to show that the $q$-generalized flexibility identity 
 for the product "$\star$" is equivalent to the fact that  the sixtuple $(\A(q), \A(q)^*, R_{\cdot}^*, L_{\cdot}^*, R_{\circ}^*, L_{\circ}^*)$ is   a
 matched pair of $\A(q)$ and $\A(q)^*.$
 For all $x, y, z\in \A(q)$, and all
 $a, b, c\in \A(q)^*$, the left and right hand sides of the  associator of the bilinear product $\star$ are given, respectively, by:
 \begin{eqnarray*}
 \{(x+a)\star(y+b)\}\star(z+c)
 &=&\{ x\cdot y+R_{\circ}^*(a)y+L_{\circ}^*(b)x+a\circ b+R_{\cdot}^*(x)b \\
 && +L_{\cdot}^*(y)a\}
 \star(z+c)
 =((x\cdot y+R_{\circ}^*(a)y \\
 && +L_{\circ}^*(b)x)\cdot z+
 R_{\circ}^*(a\circ b+R_{\cdot}^*(x)b+L_{\cdot}^*(y)a)z
 \cr 
 && +L_{\circ}^*(c)(x\cdot y+R_{\circ}^*(a)y+L_{\circ}^*(b)x)
 \cr 
 && + (a\circ b+R_{\cdot}^*(x)b+
 L_{\cdot}^*(y)a)\circ c \\
 &&
 +R_{\cdot}^*(x\cdot y+R_{\circ}^*(a)y+L_{\circ}^*(b)x)c \cr 
 &&+
 L_{\cdot}^*(z)(a\circ b 
 +R_{\cdot}^*(x)b+L_{\cdot}^*(y)a)
 \cr
 &=&(x\cdot y)\cdot z+(R_{\circ}^*(a)y)\cdot z+(L_{\circ}^*(b)x)\cdot z\cr
 &&
 +R_{\circ}^*(a\circ b)z+
 R_{\circ}^*(R_{\cdot}^*(x)b)z
 +R_{\circ}^*(L_{\cdot}^*(y)a)z\\
 && +L_{\circ}^*(c)(x\cdot y)+L_{\circ}^*(c)(R_{\circ}^*(a)y)\cr 
 && 
 +L_{\circ}^*(c)(L_{\circ}^*(b)x)\\
 && + (a\circ b)\circ c+(R_{\cdot}^*(x)b)\circ c\cr 
 && +(L_{\cdot}^*(y)a)\circ c+
 R_{\cdot}^*(x\cdot y)c+
 R_{\cdot}^*(R_{\circ}^*(a)y)c\cr 
 && +R_{\cdot}^*(L_{\circ}^*(b)x)c+
 L_{\cdot}^*(z)(a\circ b)\cr 
 &&+L_{\cdot}^*(z)(R_{\cdot}^*(x)b)
 +L_{\cdot}^*(z)(L_{\cdot}^*(y)a)
 \cr
 \{(x+a)\star(y+b)\}\star(z+c)
 &=& \{(R_{\circ}^*(a)y)\cdot z+R_{\circ}^*(L_{\cdot}^*(y)a)z\}\cr 
 &&+\{(L_{\circ}^*(b)x)\cdot z  +R_{\circ}^*(R_{\cdot}^*(x)b)z\}+(x\cdot y)\cdot z\cr 
 && +R_{\circ}^*(a\circ b)z+L_{\circ}^*(c)(x\cdot y)+L_{\circ}^*(c)(R_{\circ}^*(a)y)\cr 
 &&+ L_{\circ}^*(c)(L_{\circ}^*(b)x)
 +\{(R_{\cdot}^*(x)b)\circ c\cr 
 &&+ R_{\cdot}^*(L_{\circ}^*(b)x)c\}+
 \{(L_{\cdot}^*(y)a)\circ c \cr 
 && +R_{\cdot}^*(R_{\circ}^*(a)y)c\}+
 (a\circ b)\circ c+R_{\cdot}^*(x\cdot y)c\cr
 &&
 +L_{\cdot}^*(z)(a\circ b)+L_{\cdot}^*(z)(R_{\cdot}^*(x)b)\cr 
 &&+
 L_{\cdot}^*(z)(L_{\cdot}^*(y)a)
 \end{eqnarray*}
 and,
 \begin{eqnarray*}
 (x+a)\star\{(y+b)\star(z+c)\}
 &=&(x+a)\star\{y\cdot z+R_{\circ}^*(b)z+L_{\circ}^*(c)y\\ && 
 +b\circ c+R_{\cdot}^*(y)c+L_{\cdot}^*(z)b\}\cr
 &=&x\cdot(y\cdot z+R_{\circ}^*(b)z+L_{\circ}^*(c)y)\cr 
 &&
 +R_{\circ}^*(a)(y\cdot z+R_{\circ}^*(b)z+L_{\circ}^*(c)y)
 \cr
 &&L_{\circ}^*(b\circ c+R_{\cdot}^*(y)c+L_{\cdot}^*(z)b)x\cr 
 &&
 +a\circ (b\circ c+R_{\cdot}^*(y)c+L_{\cdot}^*(z)b)
 \cr
 &&+R_{\cdot}^*(x)(b\circ c+R_{\cdot}^*(y)c+L_{\cdot}^*(z)b)
 \cr && +L_{\cdot}^*(y\cdot z+R_{\circ}^*(b)z+L_{\circ}^*(c)y)a
 \cr
 &=&x\cdot(y\cdot z)+x\cdot(R_{\circ}^*(b)z)+x\cdot(L_{\circ}^*(c)y)\cr 
 && +R_{\circ}^*(a)(y\cdot z)+
 R_{\circ}^*(a)(R_{\circ}^*(b)z)\cr && 
 +R_{\circ}^*(a)(L_{\circ}^*(c)y)+
 L_{\circ}^*(b\circ c)x+L_{\circ}^*(R_{\cdot}^*(y)c)x
 \cr
 &&
 +L_{\circ}^*(L_{\cdot}^*(z)b)x+
  a\circ (b\circ c)+
 a\circ(R_{\cdot}^*(y)c)\cr &&
 +a\circ(L_{\cdot}^*(z)b)+
 R_{\cdot}^*(x)(b\circ c)\cr && 
 +R_{\cdot}^*(x)(R_{\cdot}^*(y)c)+
 R_{\cdot}^*(x)(L_{\cdot}^*(z)b)+
 L_{\cdot}^*(y\cdot z)a\cr
 && +L_{\cdot}^*(R_{\circ}^*(b)z)a+L_{\cdot}^*(L_{\circ}^*(c)y)a.
 \cr
 (x+a)\star\{(y+b)\star(z+c)\}
 &=&
 \{
 x\cdot(R_{\circ}^*(b)z)+L_{\circ}^*(L_{\cdot}^*(z)b)x
 \}\cr 
 && 
 +\{
 x\cdot(L_{\circ}^*(c)y)+L_{\circ}^*(R_{\cdot}^*(y)c)x
 \}\cr 
 &&
 +R_{\circ}^*(a)(y\cdot z)+R_{\circ}^*(a)(R_{\circ}^*(b)z)\cr 
 && 
 +R_{\circ}^*(a)(L_{\circ}^*(c)y)+
 L_{\circ}^*(b\circ c)x
 \cr
 && +\{a\circ(R_{\cdot}^*(y)c)+
 L_{\cdot}^*(L_{\circ}^*(c)y)a \}\cr &&
 +\{a\circ(L_{\cdot}^*(z)b)+
 L_{\cdot}^*(R_{\circ}^*(b)z)a \}   
 \cr
 &&+x\cdot(y\cdot z)+ a\circ (b\circ c)+R_{\cdot}^*(x)(b\circ c)\cr &&
 +R_{\cdot}^*(x)(R_{\cdot}^*(y)c)+ R_{\cdot}^*(x)(L_{\cdot}^*(z)b)+ L_{\cdot}^*(y\cdot z)a
 \end{eqnarray*}
 Therefore, the associator can be rewritten   as:
 \begin{eqnarray*}
 \left((x+a), (y+b ,(z+c)\right)
 &=&
 \{(x+a)\star(y+b)\}\star(z+c)\cr 
 &&
 -(x+a)\star\{(y+b)\star(z+c)\}
 \cr
 \left((x+a), (y+b ,(z+c)\right)
 &=&
 \{
 (L_{\circ}^*(b)x)\cdot z+R_{\circ}^*(R_{\cdot}^*(x)b)z
 -x\cdot(R_{\circ}^*(b)z)\cr 
 &&-L_{\circ}^*(L_{\cdot}^*(z)b)x
 \}+
  \{
 L_{\circ}^*(c)(x\cdot y) -x\cdot(L_{\circ}^*(c)y)\cr 
 && -L_{\circ}^*(R_{\cdot}^*(y)c)x
  \}+
 \{
 (R_{\circ}^*(a)y)\cdot z+R_{\circ}^*(L_{\cdot}^*(y)a)z\cr 
 &&- R_{\circ}^*(a)(y\cdot z)
 \}+
 \{
 (R_{\cdot}^*(x)b)\circ c+R_{\cdot}^*(L_{\circ}^*(b)x)c \cr 
 &&- R_{\cdot}^*(x)(b\circ c)
 \}
 \cr
 && 
 +\{
 (L_{\cdot}^*(y)a)\circ c+R_{\cdot}^*(R_{\circ}^*(a)y)c-
 a\circ(R_{\cdot}^*(y)c)\cr 
 && 
 -L_{\cdot}^*(L_{\circ}^*(c)y)a
 \}
 +\{L_{\cdot}^*(z)(a\circ b)-
 a\circ(L_{\cdot}^*(z)b)\cr 
 && 
 -L_{\cdot}^*(R_{\circ}^*(b)z)a\}
 +\{R_{\circ}^*(a\circ b)z-R_{\circ}^*(a)(R_{\circ}^*(b)z)\} \cr &+& 
 \{R_{\cdot}^*(x\cdot y)c-R_{\cdot}^*(x)(R_{\cdot}^*(y)c)\}
 \cr
 &&+\{L_{\cdot}^*(z)(R_{\cdot}^*(x)b)-
 R_{\cdot}^*(x)(L_{\cdot}^*(z)b)\}\cr 
 &&
 +\{L_{\cdot}^*(z)(L_{\cdot}^*(y)a)- L_{\cdot}^*(y\cdot z)a\}
 \cr
 &&
 +(a,b, c)+(x,y,z)+
 \{
 L_{\circ}^*(c)(R_{\circ}^*(a)y)\cr 
 && -R_{\circ}^*(a)(L_{\circ}^*(c)y)
 \}+
 \{
 L_{\circ}^*(c)(L_{\circ}^*(b)x)-
 L_{\circ}^*(b\circ c)x
 \}.
 \end{eqnarray*}
 Similarily,  we have:
 \begin{eqnarray*}
 q\left(z+c, y+b ,x+a\right)
 &=&
 q\{
 (L_{\circ}^*(b)z)\cdot x+R_{\circ}^*(R_{\cdot}^*(z)b)x
 -z\cdot(R_{\circ}^*(b)x)\cr 
 &&- L_{\circ}^*(L_{\cdot}^*(x)b)z
 \}+ 
 q\{
 (R_{\circ}^*(c)y)\cdot x+R_{\circ}^*(L_{\cdot}^*(y)c)x\cr
 && -R_{\circ}^*(c)(y\cdot x)
 \}+
 q\{
 L_{\circ}^*(a)(z\cdot y) -z\cdot(L_{\circ}^*(a)y)\cr 
 && -L_{\circ}^*(R_{\cdot}^*(y)a)z
  \}+
 q\{L_{\cdot}^*(x)(c\circ b)-
 c\circ(L_{\cdot}^*(x)b)\cr 
 && 
 -L_{\cdot}^*(R_{\circ}^*(b)x)c
 \}
 +
 q\{
 (L_{\cdot}^*(y)c)\circ a+R_{\cdot}^*(R_{\circ}^*(c)y)a\cr 
 && 
 -c\circ(R_{\cdot}^*(y)a)-
 L_{\cdot}^*(L_{\circ}^*(a)y)c
 \}+
 \{
 (R_{\cdot}^*(z)b)\circ a\cr 
 && 
 +R_{\cdot}^*(L_{\circ}^*(b)z)a-R_{\cdot}^*(z)(b\circ a)
 \}+q(c,b, a)+q(z,y,x)
 \cr
 &&
 +q\{
 R_{\cdot}^*(z\cdot y)a-R_{\cdot}^*(z)(R_{\cdot}^*(y)a)
 \}
 \cr&& 
 +q\{
 L_{\cdot}^*(x)(R_{\cdot}^*(z)b)
 -R_{\cdot}^*(z)(L_{\cdot}^*(x)b)
 \}\cr 
 && 
  +q\{
  L_{\circ}^*(a)(L_{\circ}^*(b)z)-
  L_{\circ}^*(b\circ a)z
  \}
  \cr
 &&+ q\{
  L_{\circ}^*(a)(R_{\circ}^*(c)y)-R_{\circ}^*(c)(L_{\circ}^*(a)y)
  \}
 +q\{
 R_{\circ}^*(c\circ b)x\cr 
  && -R_{\circ}^*(c)(R_{\circ}^*(b)x)
 \}
 +
 q\{
 L_{\cdot}^*(x)(L_{\cdot}^*(y)c)-L_{\cdot}^*(y\cdot x)c
 \}.
 \end{eqnarray*}
 Hence, the  equality:
 \begin{eqnarray*}
 (x+a, y+b, z+c)=q(z+c, y+b, x+a)
 \end{eqnarray*}
 is equivalent to the following:
 \begin{eqnarray*}
 \left\lbrace
 \begin{array}{ccc}
 (R_{\cdot}^*, L_{\cdot}^*, \A(q)^*) \mbox{ and }
 (R_{\circ}^*, L_{\circ}^*, \A(q)) 
 \mbox{ are bimodules of } \A(q) \mbox{ and } \A(q)^*,\\ \mbox{ respectively}, 
 (L_{\circ}^*(b)x)\cdot z+R_{\circ}^*(R_{\cdot}^*(x)b)z
 -x\cdot(R_{\circ}^*(b)z)-L_{\circ}^*(L_{\cdot}^*(z)b)x= \\ q((L_{\circ}^*(b)z)\cdot x+R_{\circ}^*(R_{\cdot}^*(z)b)x
  -z\cdot(R_{\circ}^*(b)x)-L_{\circ}^*(L_{\cdot}^*(x)b)z),
 \cr
 L_{\circ}^*(c)(x\cdot y) -x\cdot(L_{\circ}^*(c)y)-L_{\circ}^*(R_{\cdot}^*(y)c)x=\\ q((R_{\circ}^*(c)y)\cdot x+R_{\circ}^*(L_{\cdot}^*(y)c)x-R_{\circ}^*(c)(y\cdot x)),
 \cr
 (R_{\circ}^*(a)y)\cdot z+R_{\circ}^*(L_{\cdot}^*(y)a)z-R_{\circ}^*(a)(y\cdot z)=\\q(L_{\circ}^*(a)(z\cdot y) -z\cdot(L_{\circ}^*(a)y)-L_{\circ}^*(R_{\cdot}^*(y)a)z),
 \cr
 (R_{\cdot}^*(x)b)\circ c+R_{\cdot}^*(L_{\circ}^*(b)x)c-R_{\cdot}^*(x)(b\circ c)=\\q(L_{\cdot}^*(x)(c\circ b)-
 c\circ(L_{\cdot}^*(x)b)-
 L_{\cdot}^*(R_{\circ}^*(b)x)c),
 \cr
 L_{\cdot}^*(z)(a\circ b)-
 a\circ(L_{\cdot}^*(z)b)-
 L_{\cdot}^*(R_{\circ}^*(b)z)a=\\q((R_{\cdot}^*(z)b)\circ a+R_{\cdot}^*(L_{\circ}^*(b)z)a-R_{\cdot}^*(z)(b\circ a)),
 \cr
 (L_{\cdot}^*(y)a)\circ c+R_{\cdot}^*(R_{\circ}^*(a)y)c-
 a\circ(R_{\cdot}^*(y)c)-
 L_{\cdot}^*(L_{\circ}^*(c)y)a=\\q((L_{\cdot}^*(y)c)\circ a+R_{\cdot}^*(R_{\circ}^*(c)y)a
 \cr
 -c\circ(R_{\cdot}^*(y)a)-
 L_{\cdot}^*(L_{\circ}^*(a)y)c).
 \end{array} 
 \right.
 \end{eqnarray*}
 Therefore, by  Theorem~\ref{theo_dumatchedpair}, the bilinear product "$\star$" defines a $q$-generalized flexible algebra structure on the vector space 
 $\A(q)\oplus\A(q)^*$ if and only if 
 $(\A(q), \A(q)^*, R_{\cdot}^*, L_{\cdot}^*, R_{\circ}^*, L_{\circ}^*)$ is a matched pair of the $q$-generalized flexible algebras $\A(q)$ and $\A(q)^*$. $\cqfd$
 \begin{remark}
 From Proposition~\ref{prop_sum}, we conclude that both the flexible and antiflexible algebras have the same matched pairs given on $\A(q)$ and $\A(q)^*$ for $q=\pm 1$. 
 The same result extends to  associative algebras obtained for  the parameter  $q=0.$ 
 \end{remark}
 \section{Manin triple of $q$-generalized flexible algebras and bialgebras}
 
 We start with the following definitions, consistent with analogous formulation for Lie algebras \cite{yvette}:

 \begin{definition}\label{def_Manin}
 Let $(\A(q), \cdot)$ be a $q$-generalized flexible algebra.  Suppose that there is a $q$-generalized flexible algebra structure "$\circ$" on its dual space $\A(q)^*$. 
 A Manin triple of the $q$-generalized flexible algebras $\A(q)$ and $\A(q)^*$ associated to a symmetric, non-degenerate, invariant  bilinear form $\bb$ defined 
 on the vector space $\A(q)\oplus\A(q)^*$ by: 
 \begin{eqnarray}\label{bilinearform}
 \bb(x+a, y+b)=<x,b>+<y,a>,
 \end{eqnarray} for all
 $x,y \in \A(q)$ and all $a, b\in\A(q)^*,$  where the bilinear product $<,>$ is the natural pairing between the vector spaces $\A(q)$ and $\A(q)^*,$ is a triple $(\A(q)\oplus\A(q)^*,\A(q), \A(q)^*)$ such that the bilinear product "$\star$" defined  for all $x, y\in \A(q)$ and all $a, b\in \A(q)^*$ by:
 \begin{eqnarray*}
 (x+a)\star(y+b)=(x\cdot y+R_{\circ}^*(a)y+L_{\circ}^*(b)x)+(a\circ b+R_{\cdot}^*(x)b+L_{\cdot}^*(y)a)
 \end{eqnarray*}
 realizes a $q$-generalized flexible algebra structure on $\A(q)\oplus\A(q)^*.$
 \end{definition}
 \begin{definition}
 The  triple  $(\A(q), \A_1(q), \A_2(q)),$  where:
 \begin{itemize}
 \item $\A(q)$ is  a $q$-generalized flexible algebra together with a nondegenerate, invariant  and symmetric bilinear form, and  
 \item   $ \A_1(q)$ and $ \A_2(q)$ 
 are two Lagrangian sub-$q$-generalized flexible algebras of $\A(q)$  such that $\A(q)= \A_1(q)\oplus \A_2(q)$.
 \end{itemize}
 is a $q$-generalized flexible bialgebra.
 \end{definition}
 \begin{theorem}\label{theo_manin}
  Suppose there is a $q$-generalized flexible algebra structure "$\circ$" on the dual space $\A(q)^*$. 
 Then, the sixtuple $(\A(q),\A(q)^*, R_{\cdot}^*, L_{\cdot}^*, R_{\circ}^*, L_{\circ}^*)$ is a matched pair of the generalized flexible algebras $(\A(q), \cdot)$ 
 and $(\A(q)^*, \circ)$ if and only if $(\A(q)\oplus\A(q)^*, \A(q), \A(q)^*)$ is a Manin triple of the $q$-generalized flexible algebras $(\A(q), \cdot)$ and $(\A(q)^*, \circ)$.
 \end{theorem}
 {Proof}:
 
 Let $(\A(q), \cdot)$ be a $q$-generalized flexible algebra. Assume that there is a $q$-generalized flexible algebra structure "$\circ$" on its dual vector space $\A(q)^*$.
 From Theorem~\ref{prop_sum},  the vector space $\A(q)\oplus\A(q)^*$
 has a $q$-generalized flexible algebra structure here denoted by "$\star$" given,  for all $x, y \in \A(q)$ and all $a,b \in \A(q)^*,$ by:
 \begin{eqnarray*}
 (x+a)\star(y+b)=(x\cdot y+R_{\circ}^*(a)y+L_{\circ}^*(b)x)+(a\circ b+R_{\cdot}^*(x)b+L_{\cdot}^*(y)a).
 \end{eqnarray*}
 Then, the sixtuple $(\A(q), \A(q)^*, R_{\cdot}^*, L_{\cdot}^*, R_{\circ}^*, L_{\circ}^*)$ is a  matched pair of the $q$-generalized flexible algebras $(\A(q), \cdot )$ and $(\A(q)^*, \circ)$. 
 From Definition~\ref{def_Manin}, it  remains to show that the bilinear form $\bb$ defined by using the natural pairing between $\A(q)$ 
 and its dual in  \eqref{bilinearform} satisfies the relation:
 \begin{eqnarray}\label{eq_qinv}
 \bb((x+a)\star(y+b),(z+c))=\bb((x+a),(y+b)\star(z+c)).
 \end{eqnarray}
 The left hand side of the equation \eqref{eq_qinv} is given,  for all $x, y, z \in \A(q)$ and all $a, b, c \in\A(q)^*,$  by:
 \begin{eqnarray*}
 \bb((x+a)\star(y+b), (z+c))
 &=&
 \bb(x\cdot y+R_{\circ}^*(a)y+L_{\circ}^*(b)x+a\circ b+R_{\cdot}^*(x)b\\
 &&+L_{\cdot}^*(y)a, z+c)\cr
 &=&<x\cdot y+R_{\circ}^*(a)y+L_{\circ}^*(b)x, c>\cr&&
 +<z,a\circ b+R_{\cdot}^*(x)b+L_{\cdot}^*(y)a>
 \cr
 &=&<x\cdot y,c>+<R_{\circ}^*(a)y,c>\cr 
 &+& <L_{\circ}^*(b)x,c>+<z,a\circ b>\cr 
 &+&<z,R_{\cdot}^*(x)b> +<z,L_{\cdot}^*(y)a>
 =\cr 
 && <x\cdot y,c>+<y,R_{\circ}(a)c>\cr
 &+&<x,L_{\circ}(b)c>+
 <z,a\circ b>\cr 
 &+& <R_{\cdot}(x)z,b>+ <L_{\cdot}(y)z,a>
 \cr
 \bb((x+a)\star(y+b), (z+c))
 &=&<x\cdot y,c>+<y,c\circ a>+<x,b\circ c>
 \cr&+&
 <z,a\circ b>+<z\cdot x,b>+<y\cdot z,a>.
 \end{eqnarray*}
 Then, we have:
 \begin{eqnarray}\label{eq_qform1}
 \bb((x+a)\star(y+b), (z+c))
 &=&
 <x\cdot y,c>+<y,c\circ a>+<x,b\circ c>
 \cr &+&
 <z,a\circ b>+<z\cdot x,b>\\
 && +<y\cdot z,a>.\nonumber
 \end{eqnarray}
 Besides, the right hand side of the equation \eqref{eq_qinv}  can be developed as follows:
 \begin{eqnarray*}
 \bb((x+a),(y+b)\star(z+c))
 &=&\bb((x+a), y\cdot z+R_{\circ}^*(b)z+L_{\circ}^*(c)y+b\circ c\cr
 && +R_{\cdot}^*(y)c+L_{\cdot}^*(z)b   )
 \cr
 &=&<x,b\circ c>+<x,R_{\cdot}^*(y)c> \cr 
 &+& <x,L_{\cdot}^*(z)b>+
 <y\cdot z,a>\cr 
 && +<R_{\circ}^*(b)z ,a>+<L_{\circ}^*(c)y,a>
 \cr
 &=&<x,b\circ c>+<R_{\cdot}(y)x,c>
 \cr&&
 + <L_{\cdot}(z)x,b>+
 <y\cdot z,a>\cr 
 &+& <z ,R_{\circ}(b)a>+<y,L_{\circ}(c)a>
 \cr\bb((x+a),(y+b)\star(z+c))
 &=&<x,b\circ c>+<x\cdot y,c>+<z\cdot x,b>
 \cr
 &+&
 <y\cdot z,a>+<z ,a\circ b>+<y,c\circ a>
 \end{eqnarray*}
 Hence,
 \begin{eqnarray}\label{eq_qform2}
 \bb((x+a),(y+b)\star(z+c))&=&
  <x,b\circ c>+<x\cdot y,c>+<z\cdot x,b>
 \cr&&+
  <y\cdot z,a>+<z ,a\circ b>\\
  && +<y,c\circ a>.\nonumber
 \end{eqnarray}
 Therefore, from the relations \eqref{eq_qform1} and \eqref{eq_qform2}, we have the required result. 
 $\cqfd$
 \begin{theorem}
 Suppose that there is a $q$-generalized flexible algebra structure "$\circ$" on the 
  dual space $\A(q)^*$. The following propositions are equivalent:
 \begin{enumerate}
 \item $(\A(q)\oplus\A(q)^*, \A(q), \A(q)^*, \omega)$ is a Manin triple of the $q$-generalized flexible algebras $\A(q)$ and $\A(q)^*$ with the nondegenerate
 symmetric bilinear form $\omega$ difined on $\A(q)\oplus\A(q)^*$, for all $x, y\in\A(q)$ and  all $a, b\in \A(q)^*$ by:
 $\omega(x+a, y+b):=<x, b>+<y, a>$, where $<,>$ is the natural pairing between $\A(q)$ and $\A(q)^*$.
 \item The sixtuple $(\A(q), \A(q)^*, R_{\cdot}^*, L_{\cdot}^*, R_{\circ}^*, L_{\circ}^*)$ is a matched pair of the $q$-generalized flexible algebras 
 $(\A(q), \cdot)$ and $(\A(q)^*, \circ)$.
 \item $ (\A(q), \A(q)^*)$ is a $q$-generalized flexible bialgebra.
 \end{enumerate}
 \end{theorem}
 Proof:
 By considering the Theorem~\ref{theo_manin}, we deduct $(1) \Longleftrightarrow (2)$. By the definition, we also have the equivalence $(1) \Longleftrightarrow (3)$.
 $\cqfd$

 \section{Application to octonion algebra}
 \begin{definition}
 An octonion algebra $\mathcal{O}$  is an eight dimensional vector space spanned by elements $\{e_0, e_1, \cdots e_7\}$ satisfying the following relations: 
 $
 \forall i, j, k=1, \cdots 7,
 $
 \begin{eqnarray}\label{eq_def_oct} 
 e_0^2=e_0, e_ie_0=e_i=e_0e_i, e_ie_j=-\delta_{ij}e_0+c_{ijk}e_k, 
 \end{eqnarray}
 where the fully antisymmetric structure constants $c_{ijk}$ are taken to be $1$ for the
 combination of  indexes:  $$(ijk) \in \{(124), (137), (156), (235), (267), (346), (457)\},$$ 
 with the bilinear product  given in  Table 1.
 \end{definition}
 \begin{table}[h!]
 \begin{tabular}{|c|c|c|c|c|c|c|c|c|}
 \hline
 $\curvearrowright$&$e_0$&$e_1$ &$e_2$ &$e_3$ & $e_4$& $e_5$& $e_6$&$e_7$  \\ 
 \hline
 $e_0$&$e_0$&$e_1$ &$e_2$ &$e_3$ & $e_4$& $e_5$& $e_6$&$e_7$  \\ 
 \hline
 $e_1$&$e_1$&$-e_0$ &$e_4$ &$e_7$ & $-e_2$& $e_6$& $-e_5$&$-e_3$  \\ 
 \hline
 $e_2$&$e_2$&$-e_4$ &$-e_0$ &$e_5$ & $e_1$& $-e_3$& $e_7$&$-e_6$  \\ 
 \hline
 $e_3$&$e_3$&$-e_7$ &$-e_5$ &$-e_0$ & $e_6$& $e_2$& $-e_4$&$e_1$  \\ 
 \hline
 $e_4$&$e_4$&$e_2$ &$-e_1$ &$-e_6$ & $-e_0$& $e_7$& $e_3$&$-e_5$  \\ 
 \hline
 $e_5$&$e_5$&$-e_6$ &$e_3$ &$-e_2$ & $-e_7$& $-e_0$& $e_1$&$e_4$  \\ 
 \hline
 $e_6$&$e_6$&$e_5$ &$-e_7$ &$e_4$ & $-e_3$& $-e_1$& $-e_0$&$e_2$  \\ 
 \hline
 $e_7$&$e_7$&$e_3$ &$e_6$ &$-e_1$ & $e_5$& $-e_4$& $-e_2$&$-e_0$  \\ 
 \hline
 \end{tabular}
 \caption{Multiplication table of octonion algebra.}
 \end{table}
 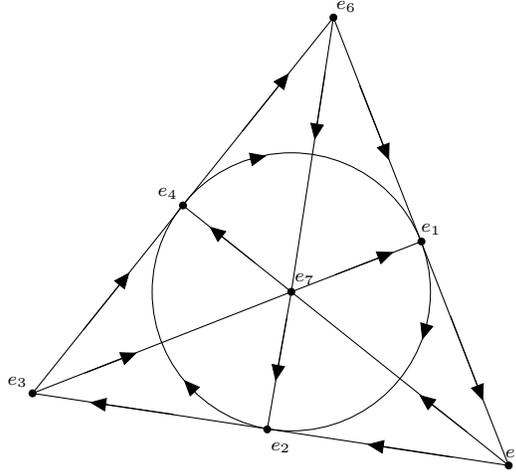
\begin{figure}
  \begin{tikzpicture}[line cap=round,line join=round,>=triangle 45,x=1.0cm,y=1.0cm]
  \draw (7,10)-- (3,5);
  \draw (3,5)-- (9.33,4.04);
  \draw (9.33,4.04)-- (7,10);
  \draw(6.44,6.35) circle (1.85cm);
  \draw (7,10)-- (6.44,6.35);
  \draw (6.44,6.35)-- (6.12,4.52);
  \draw (3,5)-- (6.44,6.35);
  \draw (6.44,6.35)-- (8.17,7.02);
  \draw (5,7.5)-- (6.44,6.35);
  \draw (6.44,6.35)-- (9.33,4.04);
  \draw [->] (3.68,5.85) -- (4.29,6.61);
  \draw [->] (5.79,8.49) -- (6.41,9.26);
  \draw [->] (6.89,9.26) -- (6.75,8.37);
  \draw [->] (6.38,6) -- (6.23,5.16);
  \draw [->] (7.35,9.09) -- (7.74,8.1);
  \draw [->] (8.56,6.02) -- (8.99,4.9);
  \draw [->] (8.58,4.15) -- (7.45,4.32);
  \draw [->] (4.87,4.72) -- (3.76,4.88);
  \draw [->] (3.7,5.27) -- (4.39,5.54);
  \draw [->] (7.03,6.58) -- (7.8,6.88);
  \draw [->] (8.77,4.49) -- (8.14,4.98);
  \draw [->] (5.95,6.74) -- (5.35,7.22);
  \draw [->] (5.89,8.11) -- (6.11,8.16);
  \draw [->] (8.24,5.91) -- (8.18,5.71);
  \draw [->] (5.16,5.02) -- (5.01,5.18);
  \begin{scriptsize}
  \fill  (7,10) circle (1.5pt);
  \draw (7.17,10.15) node {$e_6$};
  \fill  (3,5) circle (1.5pt);
  \draw (2.80,5.15) node {$e_3$};
  \fill  (9.33,4.04) circle (1.5pt);
  \draw (9.42,4.19) node {$e_5$};
  \fill  (6.44,6.35) circle (1.5pt);
  \draw (6.62,6.52) node {$e_7$};
  \fill  (5,7.5) circle (1.5pt);
  \draw (4.80,7.66) node {$e_4$};
  \fill  (8.17,7.02) circle (1.5pt);
  \draw (8.3,7.18) node {$e_1$};
  \fill  (6.12,4.52) circle (1.5pt);
  \draw (6.3,4.28) node {$e_2$};
  \end{scriptsize}
  \end{tikzpicture}
  \caption{Realization of octonion algebra}
  \end{figure}

 The associator of the octonion algebra $\mathcal{O}=Span\{e_0, e_1, \cdots e_7\}$ defined as:
 \begin{eqnarray}
   e_{ijk}:=(e_i, e_j, e_k)=(e_ie_j)e_k-e_i(e_je_k), \forall i, j, k\in \{0, 1, 2, \cdots 7\}
 \end{eqnarray}
 obeys the following relations: $\forall  i, j, k\in \{1, 2, \cdots , 7\}$,
 \begin{eqnarray}
 (e_0, e_i, e_j)=(e_i, e_0, e_j)=(e_i, e_j, e_0)=(e_i, e_i, e_j=(e_i, e_j, e_i)=\\(e_i, e_j, e_i)=0\nonumber
 \end{eqnarray}
 \begin{eqnarray}
 (e_i, e_j, e_k)= \sum_{m=1}^{7}(c_{ijm}\delta_{mk}-c_{jkm}\delta_{im})e_0+
 \sum_{n=1}^{7}\sum_{m=1}^{7}(c_{ijm}c_{mkn}-c_{jkm}c_{imn})e_n,
 \end{eqnarray}
 \begin{table}{}
  \centering
 \begin{tabular}{|c|c|c|c|c|c|c|c|c|}
 \hline
 $\curvearrowright$&$e_0$&$e_1$ &$e_2$ &$e_3$ & $e_4$& $e_5$& $e_6$&$e_7$  \\ 
 \hline
 $e_{00}$&$0$&$0$ &$0$ &$0$ & $0$& $0$& $0$&$0$  \\ 
 \hline
 $e_{01}$&$0$&$0$ &$0$ &$0$ & $0$& $0$& $0$&$0$  \\ 
 \hline
 $e_{11}$&$0$&$0$ &$0$ &$0$ & $0$& $0$& $0$&$0$   \\ 
 \hline
 $e_{12}$&$0$&$0$ &$0$ &$-2e_6$ & $0$& $2e_7$& $2e_3$&$-2e_5$  \\ 
 \hline
 $e_{22}$&$0$&$0$ &$0$ &$0$ & $0$& $0$& $0$&$0$   \\ 
 \hline
 $e_{23}$&$0$&$-2e_6$ &$0$ &$0$ & $-2e_7$& $0$& $2e_1$&$2e_4$  \\ 
 \hline
 $e_{33}$&$0$&$0$ &$0$ &$0$ & $0$& $0$& $0$&$0$  \\ 
 \hline
 $e_{34}$&$0$&$2e_5$ &$-2e_7$ &$0$ & $0$& $-2e_1$& $0$&$2e_2$  \\ 
 \hline
 $e_{44}$&$0$&$0$ &$0$ &$0$ & $0$& $0$& $0$&$0$  \\ 
 \hline
 $e_{45}$&$0$&$2e_3$ &$2e_6$ &$-2e_1$ & $0$& $0$& $-2e_2$&$0$  \\ 
 \hline
 $e_{55}$&$0$&$0$ &$0$ &$0$ & $0$& $0$& $0$&$0$ \\ 
 \hline
 $e_{56}$&$0$&$0$ &$2e_4$ &$2e_7$ & $-2e_2$& $0$& $0$&$-2e_3$  \\ 
 \hline
 $e_{66}$&$0$&$0$ &$0$ &$0$ & $0$& $0$& $0$&$0$  \\ 
 \hline
 $e_{67}$&$0$&$-2e_4$ &$0$ &$2e_5$ & $2e_1$& $-2e_3$& $0$&$0$  \\ 
 \hline
 $e_{77}$&$0$&$0$ &$0$ &$0$ & $0$& $0$& $0$&$0$  \\ 
 \hline
 \end{tabular}
 \caption{Table of composition of associator of octonion.}
 \end{table}
 where the associator is written  as $(e_i, e_j, e_k):=(e_ie_j)e_k-e_i(e_je_k):=e_{ij}e_k.$
 
 \begin{proposition}
 Let $\mathcal{O}$ be an octonion algebra with basis $\{e_0, e_1, \cdots,  e_7\}$. We have:
 \begin{enumerate}
 \item The $4$ dimensional sub-algebras, spanned by the elements $\{e_0, e_i, e_j, e_k \}$ where the index $(ijk)\in \{(124), (137), (156), (235), (267) (346), (457) \},$ are associative, i.e their associator  vanishes. So far, the associator $(e_i, e_j, e_k)$ such that indexes are repeated, or contain zero,   also vanishes, and the vector space $\{e_0, e_i, e_j, e_k\}$ does not have a sub-algebra property.
 \item Other associators $(e_i, e_j, e_k),$   where  the indexes $(ijk)$  are fully\\ skew-symmetric for the combinations :
 \\$(ijk)\in \{(123), (125),(126), (127), (234), (236), (237), (341), (342), (345),\\(347),(451), (452), (453), (456), 
 (562), (563), (564),  (567), (671), (673),\\ (674), (675) \},$   
 
  do not vanish, and are anti-left symmetric, anti-right symmetric and anti-center symmetric.
 \end{enumerate}  
 \end{proposition}
 \begin{definition}
 Let $\mathcal{O}$ be an octonion algebra. Consider the triple $(l, r, V)$, where $V$ is a finite 
 dimensional vector space, and $l, r:\mathcal{O}\rightarrow \mathfrak{gl}(V)$ are two linear
 maps. Then,  the following relations are satisfied:
 for all $e_i\in \mathcal{O}, i=0, 1, \cdots, 7,$  
 \begin{eqnarray}\label{eq_bimod_oct_1}
 l_{e_0}= \id=r_{e_0}, l_{e_i}=-r_{e_i},
 \end{eqnarray}
 \begin{eqnarray}\label{eq_bimod_oct_2}
 [r_{e_i}, l_{e_j}]=[r_{e_j}, l_{e_i}],
 \end{eqnarray}
 \begin{eqnarray}\label{eq_bimod_oct_3}
 \delta_{ij}+l_{e_i}l_{e_j}=c_{ijk}l_{e_k},
 \end{eqnarray}
 where the structure constants $c_{ijk},$ given in the equation \eqref{eq_def_oct}, are well defined.
 \end{definition}
 \begin{proposition}
 Let $\mathcal{O}$ be an octonion algebra, and $l, r: \mathcal{O}\rightarrow \mathfrak{gl}(V)$ be two linear
 maps. The couple $(l, r)$ is a bimodule of the octonion  algebra $\mathcal{O}$ if and only if there
 exists an octonion  algebra structure "$\ast$" on the semi-direct vector space $\mathcal{O}\oplus V$ given
 by 
 \begin{eqnarray*}
 (e_i + u)\ast (e_j + v) := e_ie_j + l_{e_i} v + r_{e_j} u,
 \forall e_i, e_j\in \mathcal{O},  \forall u, v \in V, i, j=0, 1, \cdots 7.
 \end{eqnarray*}
 \end{proposition}
 \begin{theorem}\label{thm_myung}
 For an octonion algebra $\mathcal{O}$ spanned by $\{e_0, e_1, \cdots, e_7\}$, the following relations are equivalent:
 \begin{enumerate}
 \item 
 \begin{eqnarray}
 [e_k, e_ie_j]=[e_k, e_i]e_j+e_i[e_k, e_j],
 \end{eqnarray}
 \item
 \begin{eqnarray}
 c_{ijm}c_{kml}=c_{kim}c_{mjl}+c_{kjm}c_{iml},
 \end{eqnarray}
 \end{enumerate}
 where the reals $C_{ijk}$ are defined in the equation\eqref{eq_def_oct},  for all $i, j, k\in \{0, 1, \cdots 7\}$. 
 \end{theorem}
  Theorem \ref{thm_myung} is known as Myung Theorem. For more details,  see \cite{myung, okubo}.
 \begin{proposition}
 Let $\mathcal{O}$ be an octonion algebra with basis $\{e_0, e_1,\cdots e_7 \}$. The following relation is satisfied: 
 \begin{eqnarray}
 2\delta_{ij}+c_{ijk}r_{e_k}+2r_{e_j}r_{e_i}
 =0,
 \end{eqnarray}
 or, equivalently,
 \begin{eqnarray} 
 2\delta_{ij}-c_{ijk}l_{e_k}+2l_{e_j}l_{e_i}
 =0,
 \end{eqnarray}
 $\forall i, j, k=0, 1, \cdots, 7$, where $c_{ijk}$ are defined in the equation\eqref{eq_def_oct}, and $l_{e_i}, r_{e_i}$ are  linear operators satisfying the relations \eqref{eq_bimod_oct_1}, \eqref{eq_bimod_oct_2} and \eqref{eq_bimod_oct_3}.
 \end{proposition}
 \section{Concluding remarks}
 In this work, we have provided a $q$-generalization of  flexible algebras,  
 including  center-symmetric, (also called antiflexible),  algebras, and their bialgebras.  
 Basic properties have been derived and discussed, as well as their connection to known algebraic
 structures investigated in the literature.  A $q$-generalization of Myung theorem  has been given.
  Main properties related to   bimodules, matched pairs, dual bimodules,  and their algebraic 
  consequences have been derived and discussed. Finally,   the equivalence between $q$-generalized 
 flexible algebras, their Manin triple and   bialgebras has been elucidated.
 
 
 

\end{document}